\documentclass[a4paper,10pt]{amsart}
\usepackage[matrix,arrow,tips,curve]{xy}
\usepackage[english]{babel}
\usepackage{amsmath}
\usepackage{amssymb}
\usepackage{pgf,tikz}
\usepackage{mathrsfs}
\usepackage{hyperref}
\usepackage{epstopdf}
\usepackage{multicol}
\usepackage{amssymb}
\usepackage{amsmath}
\usepackage{scalefnt}
\usepackage[T1]{fontenc}
\usepackage[latin1]{inputenc}
\usepackage[hyperpageref]{backref}

\input{xy}
\xyoption{all}
\newtheorem{thm}{Theorem}[section]
\newtheorem{Lemma}[thm]{Lemma}
\newtheorem{Proposition}[thm]{Proposition}
\newtheorem{Corollary}[thm]{Corollary}

\newtheorem*{thm*}{Theorem}
\newtheorem*{Proposition*}{Proposition}

\theoremstyle{definition}

\newtheorem{Definition}[thm]{Definition}
\newtheorem{Remark}[thm]{Remark}

\newtheorem{Example}[thm]{Example}

\newcommand{\Ker}{\operatorname{Ker}}

\newcommand{\Pic}{\operatorname{Pic}}

\newcommand{\codim}{\operatorname{codim}}

\newcommand{\Sing}{\operatorname{Sing}}

\begin{document}
\title[On Poincar\'e-Type Inequalities for Flags of Distributions]{On Poincar\'e-Type Inequalities for Flags of Distributions}

\author{Antonio M. Ferreira}
\address{ Antonio Marcos Ferreira da Silva \\ DMA - UFES,  Rodovia BR-101, Km 60, Bairro Litor\^aneo, S\~ao Mateus-ES, Brazil, CEP 29932-540}
\email{antonio.m.silva@ufes.br}

\author[Miguel Rodr\'iguez Peña]{Miguel Rodr\'iguez Peña}
\address{\sc Miguel Rodr\'iguez Peña\\
ICEX - UFMG\\
Departamento de Matem\'atica, Av. Ant\^onio Carlos, 6627, 31270-901  Belo Horizonte, MG, Brazil.}
\email{amrp2024@ufmg.br}

\author[Fernando Lourenço]{Fernando Lourenço}
\address{\sc Fernando Lourenço\\
DMM - UFLA\\
Campus Universitário, Lavras MG, Brazil, CEP 37200-000.}
\email{fernando.lourenco@ufla.br}

%DMM - UFLA , Campus Universitário, Lavras MG, Brazil, CEP 37200-000
%Email address: fernando.lourenco@ufla.br

%\date{\today}
%\subjclass[2010]{Primary 57R30; Secondary  14J45, 57R32, 53C12}
%\keywords{...}

\begin{abstract} In this paper, we establish enumerative formulas and derive geometric inequalities for flags of holomorphic distributions on complex manifolds. First, given a codimension-one holomorphic distribution $\mathcal{G}$ on a complex manifold $M$, we establish an upper bound for the number of independent holomorphic $p$-distributions tangent to $\mathcal{G}$ with a prescribed associated twisting line bundle 
$\mathcal{L}$. As an application, we prove Poincaré-type degree inequalities for flags on certain classes of manifolds. Finally, we investigate the Mumford--Takemoto stability of the tangent sheaves of rank-two distributions. In particular, we establish several stability results for such distributions on projective spaces, smooth weighted complete intersection threefolds, and, in particular, on smooth hypersurfaces $V\subset\mathbb{P}^4$.

\medskip
\noindent\textbf{2020 Mathematics Subject Classification}
Primary 32M25, 14J70, 14J60. \\

\noindent\textbf{Keywords} Flags; Pfaff  systems; holomorphic foliations/distributions; Poincar\'e's problem; Mumford-Takemoto stability.

\end{abstract}
\maketitle

\setcounter{tocdepth}{1}

%\tableofcontents

\section{Introduction}

A flag of holomorphic distributions on a complex manifold $M$ is a pair of distributions $\mathcal{D}_1$ and $\mathcal{D}_2$, of dimensions $k_1$ and  $k_2$ respectively, such that $k_1< k_2$ and $\mathcal{D}_1$ is a sub $\mathcal{O}_M$-module of $\mathcal{D}_2$. Moreover, if both 
$\mathcal{D}_1$ and $\mathcal{D}_2$ are integrable, then they define a flag of holomorphic foliations.

The concept of flags of foliations has recently been investigated by several authors. For instance, Feigin \cite{Feigin} studied characteristic classes of flags, Brasselet et al. \cite{bramafe} proved a Baum-Bott-type residue theorem for flags of holomorphic foliations, and Costa \cite{Cos} related the classical problem of determining when the tangent sheaf of a distribution on $\mathbb{P}^{n}$ splits to the existence of a flag. 

There are many topics which are closely related to flags and naturally appear in the theory of holomorphic foliation. For instance, we refer to a conjecture due to Marco Brunella, see \cite[p.443]{Cerv} which says that a holomorphic foliation of dimension two on $\mathbb{P}^{3}$ either admits an invariant algebraic surface, or it composes a flag of holomorphic foliations. Furthermore, we have a more general conjecture attributed to different authors, see \cite[p.2]{CerLN0}: 
\begin{center}
{\it Any codimension one holomorphic foliation $\mathcal{F}$ on $\mathbb{P}^n$, $n\geq3$, either admits a transverse projective structure with poles on some invariant hypersurface, or it is the pull-back of a holomorphic foliation on $\mathbb{P}^2$ by a rational map.}
\end{center}

\noindent In the second alternative, the leaves of $\mathcal{F}$ are foliated by the levels of the rational map.

Corr\^ea et al \cite{SoaCo} established Poincaré-type inequalities relating the degrees of holomorphic distributions and foliations forming a flag on $\mathbb{P}^n$. These inequalities are motivated by the so-called \textit{Poincar\'e problem} for foliations. 

The Poincar\'e problem originates in the work of Henri Poincar\'e \cite{Poin} and concerns the existence of a rational first integral for a holomorphic foliation $\mathcal{F}$ on $\mathbb{P}^2$. Poincaré observed that this question reduces to finding an upper bound for the degrees of the generic $\mathcal{F}$-invariant curves. Although such a bound does not exist without further assumptions, several authors have investigated the problem under suitable hypotheses, as well as various generalizations; see for instance 
\cite{BruMe,Bru1,Can,Alc,Este,Gal,Pere,So}.

In a recent article, Rodr\'iguez et al. \cite{RoFe} also considered enumerative problems for flags of distributions in the projective space 
$\mathbb{P}^n$ and, as an application, established Poincar\'e-type degree inequalities for such flags.
In the present work, we address the same type of enumerative problem in the more general setting of arbitrary varieties. Specifically, in Theorem \ref{theoA} below, we assume that $M$ is a complex manifold. 

\begin{thm}\label{theoA} Let $M$ be a complex manifold.
	Let $\mathcal{G}$ be a codimension one holomorphic distribution on $M$ induced by a section of
	$\Omega^{1}_M\otimes\mathcal{L}_1$, for some invertible sheaf $\mathcal{L}_1$ on $M$.
	Then, for each integer $p$ such that $\codim\Sing(\mathcal{G})\geq p+1$, and for any invertible sheaf 
	$\mathcal{L}$ on $M$, there are at most
	$$\wp(p,\mathcal{L})=\dim H^{0}(M, \Omega^{p-1}_M\otimes\mathcal{L}\otimes\mathcal{L}_1^{\vee}),$$
	
	\noindent independent holomorphic distributions in $H^{0}(M, \Omega^{p}_{M}\otimes\mathcal{L})$ that are tangent to $\mathcal{G}$.
\end{thm}

For example, the number $\wp(p,\mathcal{L})$ can be computed explicitly when $M$ is a smooth projective quadric hypersurface 
(see Proposition \(\ref{thq}\) and Example \(\ref{qua_n=3}\)), a Hopf manifold, or a compact toric orbifold \cite{EvM}. 
See also Remark \ref{Rk}. An application is presented in Corollary \ref{Apn}.
 
Now we consider the case of $M = \mathbb{P}^{n}$ to determine an upper bound for the number of invariant $p$-distributions using Bott's formula. In particular, this allows us to establish a Poincaré-type inequality for flags:

\begin{Corollary}\label{Co1} Let $\mathcal{G}$ be a codimension one holomorphic distribution on $\mathbb{P}^{n}$ induced by a section 
	$\omega_{1} \in H^{0}(\mathbb{P}^{n}, \Omega^{1}_{\mathbb{P}^{n}}(d_{1}))$. Then, for each integer $p$ such that  
	$\codim\Sing(\omega_{1})\geq p+1$, there are at most
	$$\wp(p,d) = \displaystyle\left\{\begin{array}{ccc}
	\binom{d-d_{1}+n-p+1}{d-d_{1}} \binom{d-d_{1}-1}{p-1},  & if & d > d_{1}+p-1  \\
	0, & if & d \leq d_{1}+p-1
	\end{array} \right.$$
	\\
	\noindent independent holomorphic distributions $\omega \in H^{0}(\mathbb{P}^{n}, \Omega^{p}_{\mathbb{P}^{n}}(d))$ 
	invariant by $\mathcal{G}$. In particular, if $\mathcal{F}\prec\mathcal{G}$ is a flag, where $\mathcal{F}$ is induced by a nontrivial section of $H^{0}(\mathbb{P}^{n},\Omega^{p}_{\mathbb{P}^{n}}(d))$, then $d_1+p\leq d$.
	Alternatively,
	$$\deg(\mathcal{G})\leq \deg(\mathcal{F})-1.$$
\end{Corollary}	

\begin{proof} For $M=\mathbb{P}^{n}$, it suffices, by Bott's formula (see Subsection \ref{Bottt}), to compute
	$$\wp(p,d)=h^{0}(\mathbb{P}^{n}, \Omega^{p-1}_{\mathbb{P}^{n}}(d-d_{1})).$$
\end{proof}

As an application, we obtain Corollary \ref{Apl}. Another type of Poincar\'e-type inequality for flags on $\mathbb{P}^{n}$ 
is given in Corollary \ref{Co2}. In the case of projective varieties, we obtain the following results.
\begin{Corollary}
Let $M$ be a smooth projective variety with $\Pic(M)\equiv\mathbb{Z}$.
Let $\mathcal{F} \prec \mathcal{G}$ be a flag of holomorphic distributions on $M$, induced by sections of
$\Omega^{p}_M(d)$ and $\Omega^{1}_M(d_1)$, respectively. Assume that  
$\codim\Sing(\mathcal{G})\geq p+1$. Then the following statements hold:
\begin{enumerate}
	\item If $M=Q_{n}\subset\mathbb{P}^{n+1}$ is a quadric hypersurface, with $n>2$, then
	$$d<d_1+p-(n+1), \,\,\,\,\mbox{or}\,\,\,\,\, d>d_1+p-1, \,\,\,\,\mbox{or}\,\,\,\,\, d=d_1, \,\,\,\,\mbox{or}\,\,\,\,\, d=d_1+2p-(n+2).$$
	\item Let $\mathbb{P}(\omega)$ be a $n$-dimensional weighted projective space, and let
	$M=V(a_1,\ldots,a_m)\subset\mathbb{P}(\omega)$ be a smooth complete intersection, where $n\geq m+3$. Then
	$$d\geq d_1+p-1.$$
\end{enumerate}
\end{Corollary}

Similarly, for Hopf manifolds, we obtain the following Poincaré-type relations for flags.

\begin{Corollary}
Let $M$ be a Hopf manifold of dimension $n\geq3$.
Let $\mathcal{F} \prec \mathcal{G}$ be a flag of holomorphic distributions on $M$, induced by sections of
$\Omega^{p}_M\otimes\mathcal{L}_a$ and $\Omega^{1}_M\otimes\mathcal{L}_b$, respectively. Assume that  
$\codim\Sing(\mathcal{G})\geq p+1$. Then the following statements hold:
\begin{enumerate}
    \item[(1)] If $M$ is classical, then $a=b\,\mu^m$, where $m\in\mathbb{N}$, $m\geq p-1$.
    \item[(2)] If $M$ is no-resonance, then
    $$
    a=b\,\mu_1^{m_1}\mu_2^{m_2}\cdots\mu_n^{m_n},
    $$
    where $m_i\in\mathbb{N}$ and there exist
    $j_1,\ldots,j_k\in\{1,\ldots,n\}$, such that
    $m_{j_1},\ldots,m_{j_k}\geq 1$.

    \item[(3)] If $M$ is weak no-resonance, then
    $$
    a=b\,\mu_1^{m_1}\mu_2^{m_2}\cdots\mu_n^{m_n},
    $$
    with $m_j\in\mathbb{N}$ for all $j=1,\ldots,n$,
    $
    m_1+m_2+\cdots+m_r=t,
    $
    and there exist
    $i_1,\ldots,i_{k-t}\geq r+1$ such that
    $
    m_{i_1}\geq 1,\ldots,m_{i_{k-t}}\geq 1.
    $
\end{enumerate}
\end{Corollary}

\vspace{1mm}

As a second application of our results, we study the Mumford--Takemoto stability of the tangent sheaves of two-dimensional distributions.
The concept of slope-stability was first introduced by D. Mumford in \cite{Mum} in the context of the construction of moduli spaces of vector bundles over a curve. This notion was later generalized to higher-dimensional varieties by Takemoto, see \cite{Ta}. 
For further references, we refer the reader to \cite{HuLe,Ko,Ta}. 

In \cite{LPT}, the authors study stability of  holomorphic foliations with trivial canonical bundle. They proved that if $\mathcal{F}$ is a semistable  codimension-one foliation on a polarized projective complex manifold $(X,H)$ with $\mathcal{K}_{\mathcal{F}}\cdot H^{n-1}=0$, then either $\mathcal{F}$ is $p$-closed or it is induced by a closed rational $1$-form with coefficients in a flat line bundle and 
and whose zero set has no divisorial components.

In \cite{OmCoJa1}, Calvo, Corrêa, and Jardim proved that every codimension-one holomorphic distribution $\mathcal{D}$ on $\mathbb{P}^{3}$ with only isolated singularities is stable. In Corollary \ref{ton}, we provide a direct proof of this result.
On the other hand, for Fano threefolds with Picard number one, the authors of \cite{ams} show that the tangent sheaf of a codimension-one locally free distribution with trivial canonical bundle is either stable or splits.

In \cite{Na} A. Napame studies the stability of the logarithmic tangent sheaf $TM(-\log D)$ associated with a log pair
$(M,D)$, where $M$ is a normal toric variety and $D$ a reduced Weil divisor. He provides a complete description of the divisors $D$ and
polarizations $\mathcal{L}$ for which $TM(-\log D)$ is (semi)stable with respect to $\mathcal{L}$ when $M$ has Picard number one or two.

Recently, Barbassa \cite{Ba} studied codimension-one distributions on $\mathbb{P}^{3}$ whose tangent sheaves are nonsplit and unstable. Moreover, he related the order of stability to the degree of the induced subfoliation by curves of degree one. 
In particular, he classified the codimension-one distributions attaining the maximal possible order of instability.

Since studying the stability of the tangent sheaf of a distribution on a projective variety is a challenging problem, we focus on rank-two distributions. We establish our main stability results in this setting. 
For projective spaces, we have the following result (see also Corollary \ref{ton} and Proposition \ref{Prop}).

\begin{thm*} Let $\mathcal{G}=\mathcal{G}_1\cap\mathcal{G}_2$ be a complete intersection distribution 
on $\mathbb{P}^4$, i.e., in homogeneous coordinates $\mathcal{G}$ is defined by 
$$\omega=\omega_1\wedge\omega_2,$$ 
where each $\mathcal{G}_i$ is a codimension one distribution induced by the homogeneous $1$-form $\omega_i$. 
Assume that each singular set of $\mathcal{G}_i$ consists only of isolated singularities. Then $\mathcal{G}$ is stable.
\end{thm*}

Finally, for smooth three-dimensional weighted complete intersections and smooth hypersurfaces $V\subset\mathbb{P}^4$, 
we obtain the following results.

\begin{thm*}
Let $\mathbb{P}(\omega)$ be a well-formed $n$-dimensional weighted projective space.
Let $V=V(a_1,\ldots,a_{n-3})\subset\mathbb{P}(\omega)$ be a smooth weighted complete intersection threefold, and let
$\mathcal{D}$ be a two-dimensional distribution on $V$. If $\mathcal{D}$ has only isolated singularities and
$\left|\omega\right|\leq \left|a\right|+2$, then $\mathcal{D}$ is stable. 
\end{thm*}
 
\begin{thm*}
Let $V\subset\mathbb{P}^4$ be a smooth hypersurface, and let
$\mathcal{D}$ be a two-dimensional distribution on $V$ with only isolated singularities.
Then $\mathcal{D}$ is stable in each of the following cases:
\begin{enumerate}
	\item $\deg(V)\neq2$;
	\item $\deg(V)=2$ and $\deg(\mathcal{D})>5$. Alternatively, $\mathcal{D}$ is semi-stable if $\deg(\mathcal{D})\geq5$.
\end{enumerate}
\end{thm*}

%\hfill

%\subsection*{Acknowledgments} ... %We are very grateful to K. Saito for the stimulating discussions we have had.

%\hfill

%\subsection*{Notation and Conventions}
%...

%%%%%%%%%%%%%%%%%%%%%%%%%%%%%%%%%%%%%%%%%%%%%%%%%%%%%%%%%%

\section{Preliminaries} \label{section:preliminaries}

\subsection{Holomorphic Distributions} 

A holomorphic singular Pfaff system $\mathcal{F}$, of codimension $k$, on a complex manifold $M$, is a non-trivial section $\omega\in H^0(M,\Omega_M ^k\otimes \mathcal{L})$, where $\mathcal{L}$ is a holomorphic line bundle on $M$, and $\Omega_{M}^{k}$ denotes the sheaf of holomorphic $k$-forms on $M$.

 %In \cite{pena} the authors studies inequalities and enumerative formulas for flags of Pfaff Systems on $\mathbb{P}^n$ and establish a bound for the number of independent twisted $r$-forms that leave invariant a one-dimensional holomorphic foliation. 

The singular set of $\mathcal{F}$ is defined by 
$$ \Sing(\mathcal{F}):=\Sing(\omega)=\{ z\in M; \, \omega(z)=0 \}$$
For more details about Pfaff systems see \cite{Covi,Coma,Este}.

A singular holomorphic distribution of codimension $k$ on $M$ is a coherent subsheaf $\mathcal{G}$ of $\Omega_{M} ^1$ 	of rank $k$. The $k$-th wedge product of the inclusion $\mathcal{G}\subset \Omega_M^1$ gives rise to a twisted holomorphic differential $k$-form with coefficients in the line bundle $\det(\mathcal{G})$. That is, a distribution of codimension $k$ induces a Pfaff system 
$$\omega\in H^0(M,\Omega_M ^k\otimes \det(\mathcal{G}))$$
of codimension $k$ which is locally decomposable outside its singular set. That is, for each point $p\in M\setminus \Sing(\omega)$ there exists a neighborhood $U$ and holomorphic $1-$forms $\omega_1, \omega_2, \ldots , \omega_k \in H^0(U,\Omega_U^1)$ such that 
$$\omega_U= \omega_1\wedge \omega_2 \wedge \cdots \wedge \omega_k$$ 
%%%

\subsubsection{Holomorphic Distributions on $\mathbb{P}^{n}$}

Let $\omega \in H^0(\mathbb{P}^n,\Omega_{\mathbb{P}^n}^k(r))$ be a holomorphic Pfaff system of codimension $k$ on $\mathbb{P}^n$. Taking a generic non-invariant linearly embedded subspace $i:H\simeq\mathbb{P}^k\hookrightarrow\mathbb{P}^n$ we have an induced non-trivial section
\[
i^*\omega\in H^0(H,\Omega_H^k(r))\simeq H^0(\mathbb{P}^k,\mathcal{O}_{\mathbb{P}^k}(-k-1+r)),
\]
since $\Omega_{\mathbb{P}^k}^k=\mathcal{O}_{\mathbb{P}^k}(-k-1)$. The tangency set between $\omega$ and $H$, denoted by $Z(i^*\omega)$, is defined as the hypersurface of zeros of $i^*\omega$ on $H$. The \emph{degree} of $\omega$, denoted by $\deg(\omega)$, is defined as the degree of $Z(i^*\omega)$ in $H$ and, therefore, is given by
\[
\deg(\omega)=-k-1+r.
\]

So, we can say that a Pfaff system of codimension $k$ and degree $d$ on $\mathbb{P}^{n}$ can be induced by a global section $\omega\in H^0(\mathbb{P}^n,\Omega_{\mathbb{P}^n}^k(d+k+1))$, where $\deg(\omega)=d$, or, as a polynomial $k$-form on $\mathbb{C}^{n+1}$ with homogeneous coefficients of degree $d+1$, see \cite{CMS}.

\subsubsection{Flag of Holomorphic Distributions}

We say that two Pfaff systems, $\mathcal{F}$ and $\mathcal{G}$ on $M$,  define a flag  if 
$\ker(w_{\mathcal{F}})\subset \Ker(w_{\mathcal{G}})$, and we denote this  by $\mathcal{F} \prec \mathcal{G},$ where 
$$\Ker(\omega) := \{ v \in TM; \, i_{v}(\omega)=0 \}. $$

Let $\mathcal{F}$ be a one-dimensional foliation induced by the vector field $X \in H^{0}(M, TM\otimes \mathcal{L}_{1})$ and let $\mathcal{G}$ be a distribution of codimension $k$ 
induced by $\omega\in H^{0}(M, \Omega_{M}^{k}\otimes \mathcal{L}_{2})$, with $1\leq k <n-1$. 
We can show that if $i_{X}\omega = 0,$ then $\mathcal{F}$ and $\mathcal{G}$ define a $2$-flag on $M.$

For a flag $\mathcal{F} \prec \mathcal{G}$ of holomorphic foliations, we have that for each point outside their singular set, the inclusion relation $T_{x}\mathcal{F} \subset T_{x}\mathcal{G}$ holds, giving that the leaves of
$\mathcal{F}$ are contained in leaves of $\mathcal{G},$ see \cite{bramafe, mol}. In this case,
we say that $\mathcal{F}$ leaves $\mathcal{G}$ invariant or that $\mathcal{F}$ is tangent to $\mathcal{G}$, or that $\mathcal{G}$ is invariant by $\mathcal{F}$. Naturally, we define the singular set of the flag $\mathcal{F}_{1}\prec\mathcal{F}_{2}$ as the following analytic set
$$\Sing(\mathcal{F}) \cup \Sing(\mathcal{G}).$$

If we have a flag of holomorphic foliations $\mathcal{F}_{1}\prec\mathcal{F}_{2}$. Thus, each foliation induces a short exact sequence of sheaves for $i=1,2.$
$$0 \longrightarrow \mathcal{F}_{i} \longrightarrow TM \longrightarrow N_{i} \longrightarrow 0,$$

\noindent with $N_{i}= TM/\mathcal{F}_{i}  $ the normal sheaf of the foliation $\mathcal{F}_{i}.$ If we denote by $N_{12}$ the relative quotient sheaf $\mathcal{F}_{2}/ \mathcal{F}_{1}$ and by $N_{\mathcal{F}}$ the normal sheaf $N_{12}\oplus N_{1}$, we have the turtle diagram

$$ \xymatrix{ 0 \ar[rd]  &  & 0 \ar[ld] &  & 0  \\
   & \mathcal{F}_{1} \ar[rd] \ar[dd] &  & N_{2} \ar[lu] \ar[ru] &  \\
   &  & TM \ar[ru] \ar[rd]  &  &  \\
   & \mathcal{F}_{2} \ar[ru]   \ar[rd] &  & N_{1} \ar[uu]  \ar[rd] &  \\
 0 \ar[ru] &  & N_{12} \ar[ru] \ar[rd]  &  & 0 \\
  & 0 \ar[ru] &   &  0 &  }$$

This is very useful to compute characteristic classes of flag, see \cite{bramafe,FeLo}.

\subsubsection{Induced distributions}

Set $n\geq m+2$.
Let $\mathbb{P}(\omega)$ be a well formed $n$-dimensional weighted projective space, 
and $V=V(a_1,\ldots,a_m)\subset\mathbb{P}(\omega)$ be a smooth weighted complete intersection.
By \cite[Theorem 3.2.4]{Do}, $\mbox{Pic}(V)\cong\mathbb{Z}$ if $n-m\geq3$. 
Furthermore $\omega_V\cong\mathcal{O}_V(\sum a_i-\sum \omega_j)$.

The inclusion map $i:V\hookrightarrow\mathbb{P}(\omega)$ induces a restriction map
$$i^{\ast}:H^{0}(\mathbb{P}(\omega),\Omega_{\mathbb{P}(\omega)}^1\otimes\mathcal{O}_{\mathbb{P}(\omega)}(d))
\rightarrow H^{0}(V,\Omega_{V}^1(d)).$$

By Bott-type vanishing formulas for $V$, $H^{0}(V,\Omega_{V}^1(d))=0$ for $d<1$; 
see Subsection \ref{Bottt}. Therefore, we assume that every distribution $\mathcal{D}$ on $V$ is induced by a non-trivial section
$\omega$ in the image of $i^{\ast}$, with $d\geq1$.

\subsection{Hopf manifolds} Let $f(z_1, z_2, \ldots, z_n) = (\mu_1 z_1, \mu_2 z_2, \ldots, \mu_n z_n)$ be a diagonal contraction in $\mathbb{C}^n$, where $0 < |\mu_i| < 1$ for all $1 \leqslant i \leqslant n$, with $n\geq2$. The quotient space $M = (\mathbb{C}^n - \{0\}) / \langle f \rangle$ is a compact complex manifold of dimension $n$ called a diagonal Hopf manifold. When $\mu_1 = \cdots = \mu_n$ we say that $M$ is a classical Hopf manifold that was first studied by Heinz Hopf in 1948, see \cite{Hopf}. Hopf showed that $M$ is diffeomorphic to the product of odd spheres $S^1 \times S^{2n-1}$ and has a complex structure which is not-K\"{a}hler. In \cite{Koda}, Kodaira classified Hopf surfaces. The geometry of Hopf manifolds has been studied by several authors, see for instance, \cite{CofeMi} and the references therein.

We consider the following types of Hopf manifolds.

\begin{Definition} We say that
\begin{itemize}

  \item[1)] $M$ is \textit{classical} if $\mu = \mu_1 = \cdots = \mu_n$.
  \item[2)] $M$ is \textit{no-resonance} if $0 < |\mu_1| \leqslant |\mu_2| \leqslant \cdots \leqslant |\mu_n| < 1$ and there is no non-trivial relation among the $\mu_i$'s of the form
  \[
  \prod_{i \in A} \mu_i^{r_i} = \prod_{j \in B} \mu_j^{r_j}, \quad r_i, r_j \in \mathbb{N}, \quad A \cap B = \varnothing, \quad A \cup B = \{1, 2, \ldots, n\}.
  \]
  \item[3)] $M$ is \textit{weak no-resonance} if $\mu_1 = \mu_2 = \cdots = \mu_r$, where $2 \leqslant r \leqslant n - 1$, and there is no non-trivial relation among the $\mu_i$'s of the form
  \[
  \prod_{i \in A} \mu_i^{r_i} = \prod_{j \in B} \mu_j^{r_j}, \quad r_i, r_j \in \mathbb{N}, \quad A \cap B = \varnothing, \quad A \cup B = \{1, r + 1, \ldots, n\}.
  \]
\end{itemize}
\end{Definition}

\noindent Note that a general Hopf manifold is always non-resonant.

A line bundle $\mathcal{L}$ on $M$ is the quotient of $(\mathbb{C}^n - \{0\}) \times \mathbb{C}$ by the operation of a representation of the fundamental group of $M$, $\varrho_{\mathcal{L}} : \pi_1(M) \simeq \mathbb{Z} \longrightarrow GL(1, \mathbb{C}) = \mathbb{C}^*$ in the following way
\[
\begin{array}{ccc}
(\mathbb{C}^n - \{0\}) \times \mathbb{C} & \longrightarrow & (\mathbb{C}^n - \{0\}) \times \mathbb{C} \\[0.4em]
(z, v) & \longmapsto & (f(z), \varrho_{\mathcal{L}}(1)\cdot v)
\end{array}
\]

We write $\mathcal{L} = \mathcal{L}_b$ for the bundle induced by the representation $\varrho_{\mathcal{L}}$ with $b = \varrho_\mathcal{L}(1)$. For $\mathcal{L}$  a holomorphic line bundle over a Hopf manifold $M$, it is knowing that $\mathcal{L}$ is a flat bundle and the Picard group $H^1(M, \mathcal{O}^*)$ is isomorphic to $\mathbb{C}^*$.

\vspace{1mm}
\subsection{Mumford-Takemoto stability notions} We say that a torsion free coherent sheaf $\mathcal{F}$ on a $n$-dimensional smooth projective variety $M$ polarized by the ample line bundle $H$ is stable if for every proper subsheaf $\mathcal{F}'$ of $\mathcal{F}$ we have $\mu(\mathcal{F}')< \mu(\mathcal{F})$, where 

$$\mu(\mathcal{F}) = \dfrac{c_1(\mathcal{F})\cdot H^{n-1}}{rk(\mathcal{F})}$$ 

\noindent it called the slope of $\mathcal{F}$ with respect to $H.$ If we have $"\leq"$ instead $"<"$ we say that $\mathcal{F}$ is semi-stable. Thus, a holomorphic distribution $\mathcal{F}$ is said to be (semi-)stable if its tangent bundle $T_{\mathcal{F}}$ is (semi-)stable.

If $\mathcal{F}$ is a stable sheaf, it then admits a filtration by slope semi-stable coherent sheaves

$$0 = \mathcal{F}_{1} \subset \mathcal{F}_{2} \subset \cdots \subset \mathcal{F}_{k}=\mathcal{F},$$

\noindent called Jordan-Holder filtration, such that its quotients have the same slope as $\mathcal{F}$, see \cite{HuyLehn, NaTi}. Of course $\mathcal{F}$ is semi-stable
if and only if $k=1$. We note that it is linked with the existence of flag of distributions.  If $\mathcal{F}$ is a distribution on $M=\mathbb{P}^{n}$ the slope definition can be read just by 

$$\mu(\mathcal{F}) = \dfrac{\dim(\mathcal{F})- \deg(\mathcal{F})}{rk(\mathcal{F})}.$$

\vspace{1mm}
%%%%%%%%%%%%%%%%%%%%%%%%%%%%%%%%%%%%%%%%%%%%%%%%%%%%%%%%%%%%%%%%%%%%%%%
\subsection{Bott's formula} \label{Bottt}

We recall below the classical Bott formulas, see \cite[p. 8]{Oko}. Let $p,q$ and $k$ be integers, with $p$ and $q$ non-negative. Then,

$$h^{q} (\mathbb{P}^{n}, \Omega_{\mathbb{P}^{n}}^{p}(k))=\left\{\begin{array}{clccc}

\binom{k+n-p}{k}\binom{k-1}{p}  & q=0,\, 0\leq p \leq n \ \ \mbox{and} \ \ k>p; \\

1 & k=0 \ \ \mbox{and} \ \ 0\leq p=q \leq n; \\

\binom{-k+p}{-k}\binom{-k-1}{n-p} & q=n,\, 0 \leq p \leq n \ \ \mbox{and} \ \ k<p-n; \\

0 & \mbox{otherwise}.
\end{array}\right.$$

In \cite[p.$174$]{Sn}, Snow established a vanishing theorem for the cohomology groups of $\Omega_M^q(k)$ on quadric hypersurfaces 
$M=Q_n\subset\mathbb{P}^{n+1}$.

\begin{thm}[Bott's formula for Quadric]
Let $M$ be a nonsingular quadric hypersurface of dimension $n$.
\begin{enumerate}
\item If $-n+q\leq k\leq q$ and $k\neq 0$ and
$k\neq -n+2q$, then $H^p(M,\Omega_M^q(k))=0$ for all $p$;

\item $H^p(M,\Omega_M^q)\neq 0$ if and only if $p=q$;

\item $H^p(M,\Omega_M^q(-n+2q))\neq 0$ if and only if
$p=n-q$;

\item If $k>q$, then $H^p(M,\Omega_M^q(k))\neq 0$
if and only if $p=0$;

\item If $k<-n+q$, then $H^p(M,\Omega_M^q(k))\neq 0$
if and only if $p=n$.
\end{enumerate}
\end{thm}

The following theorem, which describes the cohomology of smooth weighted complete intersections, 
is due to Flenner (see \cite[p.$157$]{Fle}).

\begin{thm}[Bott's formula for weighted complete intersection]
Let $M$ be a weighted complete intersection. Then,
\begin{enumerate}
    \item $h^q(M,\Omega_M^q)=1$ for $0\leq q\leq n$, $q\neq \frac{n}{2}$.

    \item $h^p(M,\Omega_M^q(t))=0$ in the following cases:
    \begin{itemize}
        \item[$a\,)$] $0<p<n$, $p+q\neq n$, and either $p\neq q$ or $t\neq 0$;
        \item[$b\,)$] $p+q>n$ and $t>q-p$;
        \item[$c\,)$] $p+q<n$ and $t<q-p$.
    \end{itemize}
\end{enumerate}
\end{thm}

\vspace{1mm}
Finally, we have the following cohomological results for Hopf manifolds.

\begin{Lemma}\cite{CofeMi} \label{Cofm}
Let $M$ be a Hopf manifold of dimension $n\geq 3$, and let $L_b$ be a line bundle
on $M$, where $b\in\mathbb{C}^*$. Then the following statements hold:
\begin{enumerate}
    \item[(1)] If $M$ is classical, then
    $
    \dim H^0(M,\Omega_M^k\otimes L_b)>0
    $
    if and only if $b=\mu^m$, where $m\in\mathbb{N}$,
    $m\geq k$.

    \item[(2)] If $M$ is no-resonance, then
    $
    \dim H^0(M,\Omega_M^k\otimes L_b)>0
    $
    if and only if
    $$
    b=\mu_1^{m_1}\mu_2^{m_2}\cdots\mu_n^{m_n},
    $$
    where $m_i\in\mathbb{N}$ and there exist
    $j_1,\ldots,j_k\in\{1,\ldots,n\}$, such that
    $m_{j_1},\ldots,m_{j_k}\geq 1$.

    \item[(3)] If $M$ is weak no-resonance, then
    $
    \dim H^0(M,\Omega_M^k\otimes L_b)>0
    $
    if and only if
    $$
    b=\mu_1^{m_1}\mu_2^{m_2}\cdots\mu_n^{m_n},
    $$
    with $m_j\in\mathbb{N}$ for all $j=1,\ldots,n$,
    $
    m_1+m_2+\cdots+m_r=t,
    $
    and there exist
    $i_1,\ldots,i_{k-t}\geq r+1$ such that
    $
    m_{i_1}\geq 1,\ldots,m_{i_{k-t}}\geq 1.
    $
\end{enumerate}
\end{Lemma}

%%%%%%%%%%%%%%%%%%%%%%%%%%%%%%%%%%%%%%%%%%%%%%%%%%%%%%%%%
%%%%%%%%%%%%%%%%%%%%%%%%%%%%%%%%%%%%%%%%%%%%%%%%%%%%%%%%%%

\vspace{1mm}

\section{Enumerative Formulas and Poincar\'e-Type Inequalities for Flags} 

In this part of the work, we consider a flag of Pfaff systems on a manifold $M$. 
First of all, we need to characterize such flags.

\begin{Lemma} \label{L} Let $M$ be a smooth projective variety.
Let $\mathcal{F}$ and $\mathcal{G}$ be holomorphic Pfaff systems on $M$, induced by $\omega \in H^{0}(M,\Omega^{p}_{M}\otimes \mathcal{O}_{M}
(D))$ and $\omega_{1} \in H^{0}(M,\Omega^{1}_{M}\otimes \mathcal{O}_{M}(D_{1}))$, respectively, for some Cartier divisors $D$ and 
$D_{1}$ on $M$, and assume that  
$\codim\Sing(\omega_{1})\geq p+1$. If $\omega \wedge \omega_{1} \equiv0$, then $\mathcal{F}$ and $\mathcal{G}$ define a flag 
$\mathcal{F} \prec \mathcal{G}$. 
Moreover, if $\mathcal{F}$ is a distribution, then the following equivalence holds:
$$\mathcal{F}\prec\mathcal{G}\,\,\Longleftrightarrow\,\,\omega \wedge \omega_1 \equiv 0.$$
\end{Lemma} 

\begin{proof} Suppose $\omega \wedge \omega_{1} \equiv 0$. Consider $\Omega_{M}^{1}\otimes \mathcal{O}_{M}(D_{1})$ that locally is a 
$\mathcal{O}_{M}$-free module of finite rank $n$, in the work of Saito in \cite{Sa}. Then we have an exact complex associated with the $1$-form 
$\omega_{1}$ up to level $p$, provided that $\codim\Sing(\omega_{1})\geq p+1$,
\begin{equation}\label{eq1}
0 \to \mathcal{O}_{M} \to  \Omega^{1}_{M}(D_{1}) \to \cdots \to \Omega^{p-1}_{M}((p-1)D_{1}) 
\stackrel{\wedge \omega_{1}\,\,\,}{\longrightarrow} 
\Omega^{p}_{M}(pD_{1}) \to \cdots \to \Omega^{n}_{M}(nD_{1})\to 0.   
\end{equation}

Tensoring the equation $(\ref{eq1})$ by $\mathcal{O}_{M}(D-pD_{1})$, we have  
\begin{equation}\label{eq11}
0 \to \mathcal{O}_{M}(D-pD_{1}) \to \cdots \to \Omega^{p-1}_{M}(D -D_{1}) 
\stackrel{\wedge \omega_{1}\,\,\,}{\longrightarrow} \Omega^{p}_{M}(D) \to \cdots \to \Omega^{n}_{M}(D+ (n-p)D_{1})\to 0.
\end{equation}

Since $\omega \wedge \omega_{1} =0,$ and the complex is exact, $\omega \in H^{0}(M, \omega_{1} \wedge \Omega_{M}^{p-1}(D-D_{1}))$. 
Therefore, $\omega = \omega_{1}\wedge \eta_{1},$ where $\eta_{1} \in H^{0}(M,\Omega^{p-1}_M(D-D_{1})).$
Then for a vector field $X \in \Ker (\omega)$, we have
$$0= i_{X}\omega = (i_{X}\omega_{1})\cdot\eta_{1} - \omega_{1}\wedge(i_{X}\eta_{1}).$$

If $g_{1} = i_{X}\omega_{1} \neq 0$, we have $\eta_{1} = \omega_{1} \wedge i_{X}\eta_{1}\,/\,g_{1}$ on $M \setminus \{g_{1} \neq 0 \}$. 
Then $\omega = 0$ on $M \setminus \{g_{1} \neq 0 \}$, an absurd.
Therefore $g_{1} = i_{X}\omega_{1} = 0$ and $\Ker (\omega) \subset \Ker (\omega_{1}),$ as desired.

On the other hand, if $\omega$ is locally decomposable (outside its singular set) and $\mathcal{F} \prec \mathcal{G},$ then $\omega \wedge \omega_{1}\equiv 0.$ In fact, if $\omega \wedge \omega_{1}\not\equiv 0,$ then $\dim (\Ker (\omega \wedge \omega_{1})) = n-p-1,$ absurd since $\Ker (\omega) \subset \Ker (\omega\wedge \omega_{1}).$

\end{proof}

\begin{thm}\label{theoB} Let $M$ be a smooth projective variety.
Let $\mathcal{G}$ be a codimension one holomorphic distribution on $M$ induced by a section $\omega_{1}$ of
$\Omega^{1}_M\otimes\mathcal{O}_M(D_{1})$, for some divisor $D_{1}$ on $M$.
Then, for each integer $p$ such that $p+1 \leq \codim\Sing(\mathcal{G})$, and for any divisor $D$ on $M$, there are at most
$$\wp(p,D)=\dim H^{0}(M, \Omega^{p-1}_M(D-D_{1})),$$

\noindent independent holomorphic distributions in $H^{0}(M, \Omega^{p}_{M}\otimes\mathcal{O}_M(D))$ tangent to $\mathcal{G}$.
\end{thm}

\begin{proof} In fact, by Lemma $\ref{L}$, a holomorphic distribution induced by $\omega \in H^{0}(M, \Omega^{p}_{M}(D))$ is tangent to 
$\mathcal{G}$ if and only if 
$\omega \wedge \omega_{1} = 0$. The condition $\omega \wedge \omega_{1} = 0$ implies that
$$\omega \in H^{0}(M,\omega_{1} \wedge \Omega^{p-1}_{M}(D - D_{1})),$$ 
as follows from the exact sequence (\ref{eq11}). Therefore,
$$\wp(p, D) = \dim H^{0}(M, \Omega^{p-1}_{M}(D - D_{1})).$$
\end{proof}

\begin{Remark}\label{Rk}
In the notation introduced above, the number of holomorphic Pfaff systems $\mathcal{F}$ tangent to $\mathcal{G}$ such that 
$\omega \wedge \omega_1 \equiv 0$ is $\wp(p, D)$.
\end{Remark}

Now, we present some basic definition in order to understand the concept of restrictions and extensions of holomorphic distribution, 
see \cite{Fi}. A codimension $q$ holomorphic distribution $D$ on smooth projective variety $M$ induced by a section
$\omega \in H^{0}(M, \Omega_{M}^{q}(D))$ is transverse to a smooth hypersurface $Y$ if the singular set of $i^{\ast}\omega$ has codimension at least two, where $i : Y \to M$ represents the natural inclusion. In this case, $i^{\ast}\omega$ determines a distribution $\mathcal{G}$ on $Y,$and we say $\mathcal{G}$ is a restriction of $\mathcal{F}$ to $Y$. If a distribution $\mathcal{G}$ over a smooth projective hypersurface $Y \subset M$ is the restriction of a distribution $\mathcal{F}$ on $M$ to $Y$, we say $\mathcal{F}$ is an extension to
$\mathcal{G}$.

Let $\mathcal{F} \prec \mathcal{G}$ be a flag of Pfaff systems on $M$.
We say that $\mathcal{F}$ \emph{extends transversally} if there exists a Pfaff system $\widetilde{\mathcal{F}}$ on $M$ such that
$\mathcal{F}=\mathcal{G}\cap\widetilde{\mathcal{F}}$. We then have the following corollary.

\begin{Corollary} \label{Apn}
Let $M$ be a smooth projective variety.
Let $\mathcal{F} \prec \mathcal{G}$ be a flag of distributions on $M$, induced by sections of
$\Omega^{p}_M\otimes\mathcal{O}_M(D)$ and $\Omega^{1}_M\otimes\mathcal{O}_M(D_{1})$, respectively, and assume that  
$\codim\Sing(\mathcal{G})\geq p+1$. Then, $\mathcal{F}$ extends, 
and the number of independent transverse extensions does not exceed $$\wp(p,D)=\dim H^{0}(M, \Omega^{p-1}_{M}(D-D_{1})).$$
\end{Corollary}	

\begin{proof}
Let $\omega \in H^{0}(M,\Omega^{p}_{M}\otimes \mathcal{O}_{M}(D))$ and 
$\omega_{1} \in H^{0}(M,\Omega^{1}_{M}\otimes \mathcal{O}_{M}(D_{1}))$ be sections defining
$\mathcal{F}$ and $\mathcal{G}$, respectively. By the proof of Lemma \ref{L}, we have $\omega\in \Ker (\wedge\,\omega_1)$ and hence
$\omega = \omega_{1}\wedge \eta_{1}$, where $\eta_{1} \in H^{0}(M,\Omega^{p-1}_M(D-D_{1}))$. 
Therefore, the corollary follows from the fact that the number of independent transverse extensions $\eta_{1}$ does not exceed $\wp(p,D)$.
\end{proof}

\vspace{0.1cm}

Next, we consider Poincar\'e-type inequalities for flags on $\mathbb{P}^{n}$.

\begin{Corollary}\label{Co1} Let $\mathcal{G}$ be a codimension one holomorphic distribution on $\mathbb{P}^{n}$ induced by a section 
$\omega_{1} \in H^{0}(\mathbb{P}^{n}, \Omega^{1}_{\mathbb{P}^{n}}(d_{1}))$. Then, for each integer $p$ such that  
$\codim\Sing(\omega_{1})\geq p+1$, there are at most
$$\wp(p,d) = \displaystyle\left\{\begin{array}{ccc}
\binom{d-d_{1}+n-p+1}{d-d_{1}} \binom{d-d_{1}-1}{p-1},  & if & d > d_{1}+p-1  \\
0, & if & d \leq d_{1}+p-1
\end{array} \right.$$
\\
\noindent independent holomorphic distributions $\omega \in H^{0}(\mathbb{P}^{n}, \Omega^{p}_{\mathbb{P}^{n}}(d))$ invariant by $\mathcal{G}$.
In particular, if $\mathcal{F}\prec\mathcal{G}$ is a flag, where $\mathcal{F}$ is induced by a nontrivial section of $H^{0}(\mathbb{P}^{n},\Omega^{p}_{\mathbb{P}^{n}}(d))$, then $d_1+p\leq d$.
Alternatively,
$$\deg(\mathcal{G})\leq \deg(\mathcal{F})-1.$$
\end{Corollary}	
	
\begin{Remark}\label{rmk}	
When $p=n-1$, the singular set of $\omega_1$ consists only of isolated singularities, and moreover
$$\wp(n-1,d) = \displaystyle\left\{\begin{array}{ccc}
\binom{d-d_{1}+2}{d-d_{1}} \binom{d-d_{1}-1}{n-2},  & if & d > d_{1}+n-2  \\
0, & if & d \leq d_{1}+n-2.
\end{array} \right.$$

Consequently, we generalize and improve \cite[Theorem 1.5]{RoFe} by removing the hypotheses stated there. 
In particular, this also strengthens \cite[Corollary 1.7]{RoFe} as follows:
\end{Remark}

\begin{Corollary} \label{Apl} There does not exist flag $\mathcal{F}\prec \mathcal{G}$ on $\mathbb{P}^{n}$ satisfying the conditions
\begin{enumerate}
 \item[a)] $\dim(\mathcal{F}) = \codim(\mathcal{G}) = 1,$
 \item[b)] both $\mathcal{F}$ and $\mathcal{G}$ have only isolated singularities.
\end{enumerate}
\begin{proof} See Corollary 1.2 in \cite{RoFe} and Remark \ref{rmk}
\end{proof}
\end{Corollary}

\vspace{0.1cm}

\begin{Corollary} \label{Co2} Let $\mathcal{F}\prec\mathcal{G}$ be a flag of holomorphic distributions on $\mathbb{P}^{n}$, 
where $\mathcal{F}$ is induced by a nontrivial section of $H^{0}(\mathbb{P}^{n},\Omega^{p}_{\mathbb{P}^{n}}(d))$, and
$\mathcal{G}=\mathcal{G}_1\cap\cdots\cap\mathcal{G}_k$ is a complete intersection distribution, i.e.,
in homogeneous coordinates $\mathcal{G}$ is defined by 
$$\omega=\omega_1\wedge\dots\wedge\omega_k,$$ 
where each $\mathcal{G}_i$ is a codimension one distribution induced by the homogeneous $1$-form $\omega_i$. 
Assume that $\mathrm{codim}\big(\mathrm{Sing}(\mathcal{G}_i)\big)\geq p+1$, for all $i=1,\ldots,k$. Then
\begin{enumerate}
	\item $\max\hspace{-0,5mm}\left\{\deg(\mathcal{G}_i)\right\}\leq\deg(\mathcal{F})-1$. 
	\vspace{1mm}	
	\item $\deg(\mathcal{G})\leq \mathrm{codim}(\mathcal{G})\cdot\deg(\mathcal{F})-1$.
\end{enumerate}
\end{Corollary} 

\begin{proof} 
Item $(1)$ follows from Corollary \ref{Co1}, since $\mathcal{F}\prec\mathcal{G}_i$ for all $i=1,\ldots,k$.
On the other hand, because $\deg(\mathcal{G})=\sum\deg(\mathcal{G}_i)+k-1$, assertion $(2)$ follows from $(1)$.
\end{proof}

Decomposable forms and complete intersection of foliations were well studied in \cite{CerLN}. \\

Now, let us consider the Poincar\'e-type inequality for flags on some smooth projective varieties.

\begin{Corollary} \label{Exam}
Let $M$ be a smooth projective variety with $\Pic(M)\equiv\mathbb{Z}$.
Let $\mathcal{F} \prec \mathcal{G}$ be a flag of holomorphic distributions on $M$, induced by sections of
$\Omega^{p}_M(d)$ and $\Omega^{1}_M(d_1)$, respectively. Assume that  
$\codim\Sing(\mathcal{G})\geq p+1$. Then the following statements hold:
\begin{enumerate}
	\item If $M=Q_{n}\subset\mathbb{P}^{n+1}$ is a quadric hypersurface, with $n>2$, then
	$$d<d_1+p-(n+1), \,\,\,\,\mbox{or}\,\,\,\,\, d>d_1+p-1, \,\,\,\,\mbox{or}\,\,\,\,\, d=d_1, \,\,\,\,\mbox{or}\,\,\,\,\, d=d_1+2p-(n+2).$$
	\item Let $\mathbb{P}(\omega)$ be a $n$-dimensional weighted projective space, and let
	$M=V(a_1,\ldots,a_m)\subset\mathbb{P}(\omega)$ be a smooth complete intersection, where $n\geq m+3$. Then
	$$d\geq d_1+p-1.$$
	%\item Let $M=\operatorname{Gr}(k,n)\neq \operatorname{Gr}(1,3)$ be the Grassmannian of $k$-planes in $\mathbb{P}^n$, where $k\geq 1$. 
	%Assume that $n\geq 2k+1$ and $p+k>\dim M=(k+1)(n-k)$. Then 
	%$$d< d_1+1 \,\,\,\,\mbox{or}\,\,\,\, d>d_1+n.$$
	%Note that the hypothesis $n\geq 2k+1$ is natural, since $\operatorname{Gr}(k,n)\cong \operatorname{Gr}(n-k-1,n)$.
\end{enumerate}
\end{Corollary}

\begin{proof} 
The Poincaré-type inequalities can be derived once the conditions on 
$d-d_1$ ensuring that $\wp(p,d)=\dim H^{0}(M, \Omega^{p-1}_{M}(d-d_1))\neq0$ are satisfied, see Subsection \ref{Bottt}.
\end{proof}

\bigskip

It is straightforward to verify that analogous results hold for complex manifolds, not necessarily projective, for instance, 

\begin{thm} Let $M$ be a complex manifold.
Let $\mathcal{G}$ be a codimension one holomorphic distribution on $M$ induced by a section of
$\Omega^{1}_M\otimes\mathcal{L}_1$, for some invertible sheaf $\mathcal{L}_1$ on $M$.
Then, for each integer $p$ such that $\codim\Sing(\mathcal{G})\geq p+1$, and for any invertible sheaf 
$\mathcal{L}$ on $M$, there are at most
$$\wp(p,\mathcal{L})=\dim H^{0}(M, \Omega^{p-1}_M\otimes\mathcal{L}\otimes\mathcal{L}_1^{\vee}),$$

\noindent independent holomorphic distributions in $H^{0}(M, \Omega^{p}_{M}\otimes\mathcal{L})$ that are tangent to $\mathcal{G}$.
\end{thm}

In the following corollary, we derive Poincaré-type relations for the following classes of
Hopf manifolds.

\begin{Corollary} Let $M$ be a Hopf manifold of dimension $n\geq3$.
Let $\mathcal{F} \prec \mathcal{G}$ be a flag of holomorphic distributions on $M$, induced by sections of
$\Omega^{p}_M\otimes\mathcal{L}_a$ and $\Omega^{1}_M\otimes\mathcal{L}_b$, respectively. Assume that  
$\codim\Sing(\mathcal{G})\geq p+1$. Then the following statements hold:
\begin{enumerate}
    \item[(1)] If $M$ is classical, then $a=b\,\mu^m$, where $m\in\mathbb{N}$, $m\geq p-1$.
    \item[(2)] If $M$ is no-resonance, then
    $$
    a=b\,\mu_1^{m_1}\mu_2^{m_2}\cdots\mu_n^{m_n},
    $$
    where $m_i\in\mathbb{N}$ and there exist
    $j_1,\ldots,j_k\in\{1,\ldots,n\}$, such that
    $m_{j_1},\ldots,m_{j_k}\geq 1$.

    \item[(3)] If $M$ is weak no-resonance, then
    $$
    a=b\,\mu_1^{m_1}\mu_2^{m_2}\cdots\mu_n^{m_n},
    $$
    with $m_j\in\mathbb{N}$ for all $j=1,\ldots,n$,
    $
    m_1+m_2+\cdots+m_r=t,
    $
    and there exist
    $i_1,\ldots,i_{k-t}\geq r+1$ such that
    $
    m_{i_1}\geq 1,\ldots,m_{i_{k-t}}\geq 1.
    $
\end{enumerate}
\end{Corollary}

\begin{proof}
The Poincaré-type relations can be derived once the conditions on 
$ab^{-1}$ ensuring that $\wp(p,a)=\dim H^{0}(M, \Omega^{p-1}_{M}\otimes\mathcal{L}_{ab^{-1}})\neq0$ are satisfied, 
see Lemma \ref{Cofm} in Subsection \ref{Bottt}.
\end{proof}

\vspace{1mm}

Now, as we specify the above number to distribution on smooth quadric hypersurface, we recall it notion, and present some properties that will be necessary throughout the text. We denote by $Q_n$ the smooth quadric hypersurface of dimension $n$ in the projective space $\mathbb{P}^{n+1}$, and recall that this variety can be considered as a compact homogeneous manifold of the form 

$$G/P,$$

\noindent where $G = S(m+2, \mathbb{C}),$ and $P$ is the maximal parabolic subgroup, see \cite{Sn}. Then, we have

\begin{Proposition}\label{thq} For $M=Q_{n}$, with $n>2$, $2\leq l \leq n$, we have

%\begin{eqnarray*}
%&&h^{0}(Q_{n},\Omega_{Q_{n}}^{l-1}(s))=\footnotesize\left\{\begin{array}{ccccc}
%\displaystyle\sum_{r=0}^{l-1}(-1)^{r}\left\{h^{0}(\mathbb{P}^{n+1}, \Omega_{\mathbb{P}^{n+1}}^{l-r-1}(s-2r)) - 
%h^{0}(\mathbb{P}^{n+1}, \Omega_{\mathbb{P}^{n+1}}^{l-r-1}(s-2r-2))\right\}, & if & s>l \\ \\
%\binom{n+2}{l}, & if & s=l  \\ \\
%0, & if & s\leq l-1.\end{array}\right.
%\end{eqnarray*}

%%%%%%%%%%%%%%%%%%%%%%%%%%%%%%%%%%%%%%%%%%
%\begin{eqnarray*}
%&&h^{0}(Q_{n},\Omega_{Q_{n}}^{l-1}(s))= \footnotesize\left\{\begin{array}{ccccc}
%\displaystyle\sum_{r=0}^{l-1}(-1)^{r}\left\{h^{0}(\mathbb{P}^{n+1}, \Omega_{\mathbb{P}^{n+1}}^{l-r-1}(s-2r)) - 
%h^{0}(\mathbb{P}^{n+1}, \Omega_{\mathbb{P}^{n+1}}^{l-r-1}(s-2r-2))\right\}, & if & s>l-1 \\ \\
%0, & if & s\leq l-1.\end{array}\right.
%\end{eqnarray*}

\begin{eqnarray*}
&&h^{0}(Q_{n},\Omega_{Q_{n}}^{l-1}(s))= \footnotesize\left\{\begin{array}{ccccc}
\displaystyle\sum_{r=0}^{l-1}(-1)^{r}\left\{h^{0}(\mathbb{P}^{n+1}, \Omega_{\mathbb{P}^{n+1}}^{l-1-r}(s-2r)) - 
h^{0}(\mathbb{P}^{n+1}, \Omega_{\mathbb{P}^{n+1}}^{l-1-r}(s-2-2r))\right\}, & if & s>l-1 \\ \\
0, & if & s\leq l-1.\end{array}\right.
\end{eqnarray*}
%%%%%%%%%%%%%%%%%%%%%%%%%%%%%%%%%%%%%%%%%%

%$$\footnotesize h^{0}(Q_{n},\Omega_{Q_{n}}^{l-1}(s))=\left\{\begin{array}{ccccc}
%0 & if & s\leq l-1 \\ \\
%\displaystyle\sum_{r=0}^{l-1}(-1)^{r}\Big [h^{0}(\mathbb{P}^{n+1}, \Omega_{\mathbb{P}^{n+1}}^{l-r-1}(s-2r)) - h^{0}(\mathbb{P}^{n+1}, \Omega_{\mathbb{P}^{n+1}}^{l-r-1}(s-2r-2))\Big ] & if & \begin{array}{ccc} s>l-1, \\ (s,l) \neq (2,2)\end{array}  \\ \\
%h^{0}(\mathbb{P}^{n+1}, \Omega_{\mathbb{P}^{n+1}}^{1}(2)) & if & s=2 \ and \ l=2
%\end{array} \right.$$
%\noindent for $2\leq l \leq n$ and $s>l-1$. % Moreover, if $s \leq l-1$ this number is zero.
\end{Proposition}

\begin{proof} From Bott's formula for quadric, see Subsection \ref{Bottt}, we have that $h^{0}(Q_{n},\Omega_{Q_{n}}^{l-1}(s)) =0$ for $2\leq l \leq n,$ and $s\leq l-1$. We will consider now the case $s>l-1.$ Consider the short exact sequences below, % see \cite{Ca}.

\begin{eqnarray}\label{seq1}
0\to \Omega_{\mathbb{P}^{n+1}}^{q}(k-2) \to \Omega_{\mathbb{P}^{n+1}}^{q}(k)\to \Omega_{\mathbb{P}^{n+1}}^{q}(k)|_{Q_{n}} \to 0
\end{eqnarray}

\noindent and

\begin{eqnarray}\label{seq2}
0\to \Omega_{Q_{n}}^{q-1}(k-2) \to \Omega_{\mathbb{P}^{n+1}}^{q}(k)|_{Q_{n}}\to \Omega_{Q_{n}}^{q}(k) \to 0.
\end{eqnarray}

For the case $l=2$, we take cohomology groups in the sequence (\ref{seq2}) with $q =1$ and $k=s $ Then

\begin{eqnarray}\label{eq3}
h^{0}(Q_{n}, \Omega^{1}_{Q_{n}}(s)) = h^{0}(Q_{n}, \Omega^{1}_{\mathbb{P}^{n+1}}(s)|_{Q_{n}}) - h^{0}(Q_{n}, \mathcal{O}_{Q_{n}}(s-2)),
\end{eqnarray}

\noindent since $h^{1}(Q_{n}, \mathcal{O}_{Q_{n}}(s-2)) = 0,$ by Bott's formula for quadric.

Now, consider $q=1$ and $k=s$ in the sequence $(\ref{seq1})$ and taking the cohomology groups

\begin{eqnarray}\label{eq5}
h^{0}(Q_{n}, \Omega^{1}_{\mathbb{P}^{n+1}}(s)|_{Q_{n}}) = h^{0}(\mathbb{P}^{n+1}, \Omega^{1}_{\mathbb{P}^{n+1}}(s)) - h^{0}(\mathbb{P}^{n+1}, \Omega_{\mathbb{P}^{n+1}}^{1}(s-2)),
\end{eqnarray}

\noindent since $h^{1}(\mathbb{P}^{n+1}, \Omega_{\mathbb{P}^{n+1}}^{1}(s-2)) = 0$ for $s \neq 2.$

\noindent Then for $s=2,$ we get

\begin{eqnarray}\label{eq6}
h^{0}(Q_{n}, \Omega^{1}_{\mathbb{P}^{n+1}}(2)|_{Q_{n}}) = h^{0}(\mathbb{P}^{n+1}, \Omega^{1}_{\mathbb{P}^{n+1}}(2)) + h^{1}(\mathbb{P}^{n+1}, \Omega_{\mathbb{P}^{n+1}}^{1}).
\end{eqnarray}

\noindent Consider $q=0$ and $k=s-2$ in the sequence $(\ref{seq1})$ and observe that 
$\mathcal{O}_{\mathbb{P}^{n+1}}(s-2)|_{Q_{n}} = \mathcal{O}_{Q_{n}}(s-2).$ We have

\begin{eqnarray}\label{eq7}
h^{0}(Q_{n}, \mathcal{O}_{Q_{n}}(s-2))= h^{0}(\mathbb{P}^{n+1}, \mathcal{O}_{\mathbb{P}^{n+1}}(s-2)) - h^{0}(\mathbb{P}^{n+1}, \mathcal{O}_{\mathbb{P}^{n+1}}(s-4)).
\end{eqnarray}

%$$h^{0}(Q_{n}, \mathcal{O}_{Q_{n}}(s-2))= \left\{\begin{array}{ccc}
%h^{0}(\mathbb{P}^{n+1}, \mathcal{O}_{\mathbb{P}^{n+1}}(s-2)) - h^{0}(\mathbb{P}^{n+1}, \mathcal{O}_{\mathbb{P}^{n+1}}(s-4)) , & if & s-4\geq0 \\ \\
%
%
%\end{array} \right.
%$$

\noindent Substituting the expressions $(\ref{eq5})$, $(\ref{eq6})$ and $(\ref{eq7})$ in $(\ref{eq3})$ we get the result for $l=2.$
Now, the result follows by induction on $l$.

\end{proof}

\begin{Example}\label{qua_n=3} For the specific case  $n=3$ and $p=2$, we have

$$\wp(2,d) =  \left\{\begin{array}{ccccc}
\frac{1}{2}(d-d_1-1)(d-d_1+2)(2d-2d_1+1), & if & d> d_1+4 \\[0.2cm]
90, & if & d=d_1+4 \\[0.2cm]
35, & if & d=d_1+3 \\[0.2cm]
9, & if & d=d_1+2 \\[0.2cm]
0, & if & d\leq d_1+1.\end{array}\right.$$
In fact, for $d-d_1 > 4$, all four terms of $\wp(2,d)$ are nonzero, and
$$\wp(2,d) = \binom{d-d_1+3}{d-d_1}(d-d_1-1) -\binom{d-d_1+1}{d-d_1-2}(d-d_1-3) -\binom{d-d_1+2}{d-d_1-2} + \binom{d-d_1}{d-d_1-4}.$$
\end{Example} 

\vspace{2mm}

\section{Stability of Rank-Two Holomorphic Distributions}

Now, as an application, we consider the stability problem for holomorphic distributions on some projective manifolds.
In the next Corollary, we provide a direct proof that recovers the result from \cite[p.9014]{OmCoJa1}.

\begin{Corollary} \label{ton} Let $\mathcal{G}$ be a codimension one holomorphic distribution on $\mathbb{P}^{3}$. Assume that
$\mathrm{Sing}(\mathcal{G})$ has only isolated singularities. Then $\mathcal{G}$ is stable. 
\end{Corollary}

\begin{proof}
The slope of $\mathcal{G}$ is $\mu(T_{\mathcal{G}})=\big(2-\deg(\mathcal{G})\big)/2$.
Given $\mathcal{F}$ any nontrivial subsheaf of $\mathcal{G}$, its slope is $\mu(T_{\mathcal{F}})=1-\deg(\mathcal{F})$.
By Corollary $\ref{Co1}$, we have 
$$\deg(\mathcal{G})\leq \deg(\mathcal{F})-1,$$
which implies $\mu(T_{\mathcal{F}})<\mu(T_{\mathcal{G}})$.
\end{proof}

Regarding the stability of distributions on projective spaces, we also have the following results.

\begin{thm} Let $\mathcal{G}=\mathcal{G}_1\cap\mathcal{G}_2$ be a complete intersection distribution 
on $\mathbb{P}^4$, i.e., in homogeneous coordinates $\mathcal{G}$ is defined by 
$$\omega=\omega_1\wedge\omega_2,$$ 
where each $\mathcal{G}_i$ is a codimension one distribution induced by the homogeneous $1$-form $\omega_i$. 
Assume that each singular set of $\mathcal{G}_i$ consists only of isolated singularities. Then $\mathcal{G}$ is stable.
\end{thm}

\begin{proof} 
The slope of $\mathcal{G}$ is $\mu(T_{\mathcal{G}})=\big(2-\deg(\mathcal{G})\big)/2$. 
For any nontrivial subsheaf $\mathcal{F}$ of $\mathcal{G}$, its slope is $\mu(T_{\mathcal{F}})=1-\deg(\mathcal{F})$.
By Corollary $\ref{Co2}$, we have 
$$\deg(\mathcal{G})\leq 2\deg(\mathcal{F})-1,$$ 
which implies $\mu(T_{\mathcal{F}})<\mu(T_{\mathcal{G}})$.
\end{proof}

\begin{Proposition} \label{Prop}
Let $\mathcal{G}$ be a rank-two holomorphic distribution on $\mathbb{P}^n$, $n\geq 3$. 
Assume that $\Sing(\mathcal{G})$ has a Baum-Kupka component. Then $\mathcal{G}$ is semi-stable. 
\end{Proposition}

\begin{proof} 
The slope of $\mathcal{G}$ is $\mu(T_{\mathcal{G}})=\big(2-\deg(\mathcal{G})\big)/2$. 
For any nontrivial subsheaf $\mathcal{F}$ of $\mathcal{G}$, its slope is $\mu(T_{\mathcal{F}})=1-\deg(\mathcal{F})$.
By \cite[Theorem 1.3]{SoaCo}, we have 
$$\deg(\mathcal{G})\leq \deg(\mathcal{F}),$$ 
which implies $\mu(T_{\mathcal{F}})\leq\mu(T_{\mathcal{G}})$.
\end{proof}

Finally, we establish the following stability results for two-dimensional distributions on smooth weighted complete intersection threefolds and, in particular, on smooth hypersurfaces $V\subset\mathbb{P}^4$.

\begin{thm} \label{meme}
Let $\mathbb{P}(\omega)$ be a well-formed $n$-dimensional weighted projective space.
Let $V=V(a_1,\ldots,a_{n-3})\subset\mathbb{P}(\omega)$, be a smooth weighted complete intersection threefold, and let
$\mathcal{D}$ be a two-dimensional distribution on $V$. If $\mathcal{D}$ has only isolated singularities and
$\left|\omega\right|\leq \left|a\right|+2$, then $\mathcal{D}$ is stable. 
\end{thm}

\begin{proof} 
Similarly, as in the case of a smooth complete intersection in complex projective space,
the total Chern class of $V$ is
$$
c(V)=\frac{c\hspace{-0,5mm}\left(\mathcal{T}\mathbb{P}(\omega)|_V\right)}
{c\hspace{-0,5mm}\left(\mathcal{N}_{\mathbb{P}(\omega)/V}\right)}
=\frac{c\hspace{-0,5mm}\left(\bigoplus_{i=0}^n\mathcal{O}_{\mathbb{P}(\omega)}(\omega_i)|_V\right)}
{c\hspace{-0,5mm}\left(\bigoplus_{i=1}^{n-3}\mathcal{O}_{\mathbb{P}(\omega)}(a_i)|_V\right)}.
$$
Applying Whitney's formula, we get 
$$c(V)=\frac{\prod_{i=0}^n\left(1+\omega_i\,\eta\right)}{\prod_{i=0}^{n-3}\left(1+a_i\,\eta\right)},$$
\noindent where $\eta=c_1\hspace{-1mm}\left(i^{\ast}\mathcal{O}_{\mathbb{P}(\omega)}(1)\right)$, and $i:V\hookrightarrow\mathbb{P}(\omega)$ 
is the inclusion. Then 
\begin{eqnarray*}
c(V)&=&\left\{\sum_{j=0}^{\infty}(-1)^{j}\mathcal{W}_j(a)\,\eta^{\,j}\right\}
\left\{\sum_{k=0}^{n+1}\texttt{C}_{k}(\omega)\,\eta^k\right\} \\
&=& \sum_{i=0}^{3}\left\{\sum_{j+k=i}(-1)^{j}\mathcal{W}_j(a)\texttt{C}_{k}(\omega)\right\}\eta^{\,i} \\
&=& \sum_{i=0}^{3}\left\{\sum_{j=0}^{i}(-1)^{j}\mathcal{W}_j(a)\texttt{C}_{i-j}(\omega)\right\}\eta^{\,i},
\end{eqnarray*}
so that
\begin{eqnarray} \label{eq1st}
c_1(V)=\left\{\sum_{j=0}^{1}(-1)^{j}\texttt{C}_{1-j}(\omega)\mathcal{W}_j(a)\right\}\eta
=\left\{\sum_{i=0}^n \omega_i-\sum_{i=1}^{n-3} a_i\right\}\eta.
\end{eqnarray}

\vspace{0.1cm}

On the other hand, given a two-dimensional distribution $\mathcal{D}$ on $V$, we can consider the global section
$$w\in H^0\big(V,\Omega^1_{V}(d_1)\big)$$
that induces $\mathcal{D}$.
Since the singular set of $\mathcal{D}$ has codimension at least two we obtain the adjunction formula 
$$
K_V=K_{\mathcal{D}}\otimes \det N^{\ast}_{\mathcal{D}},
$$
where $K_{\mathcal{D}}=\left(\wedge^2\,T^{\ast}_{\mathcal{D}}\right)^{\ast\ast}$ is the canonical bundle and $N_{\mathcal{D}}$ the normal bundle of $\mathcal{D}$. Note that 
$\det(N_{\mathcal{D}})=\mathcal{O}_{V}(d_1)$.
In particular, by (\ref{eq1st}), we have 
$$c_1(T_{\mathcal{D}})=c_1\big(\mathcal{O}_{V}(\left|\omega\right|-\left|a\right|-d_1)\big).$$
Therefore (see \cite[Section 4]{Ro}),
$$\mu(T_{\mathcal{D}})=\frac{c_1(T_{\mathcal{D}})\cdot \mathcal{O}_{V}(1)^{2}}{rk(T_{\mathcal{D}})}=
\frac{(\left|\omega\right|-\left|a\right|-d_1)\,a_0\cdots a_{n-3}}{2\,\omega_0\cdots \omega_n}.$$

Now, let us consider a global section 
$X\in H^0\big(V,\Omega^{2}_{V}(d)\big)$ 
tangent to $\mathcal{D}$. In particular, $d\geq2$, since $H^0\big(V,\Omega^2_{V}(d)\big)\neq0$; see Subsection \ref{Bottt}.
Similarly, we have 
$$\mu(T_{X})=\frac{(\left|\omega\right|-\left|a\right|-d)\,a_0\cdots a_{n-3}}{\omega_0\cdots \omega_n}.$$
Then, by Corollary \ref{Exam} $(2)$, we have $d-d_1\geq1$. Moreover, by the hypothesis $\left|\omega\right|\leq \left|a\right|+2$, we obtain 
$2d>\left|\omega\right|-\left|a\right|+d_1$. Therefore, $\mu(T_{X})<\mu(T_{\mathcal{D}})$.
\end{proof}

\begin{Definition}
Let $V\neq \mathbb{P}^n$ be a smooth $n$-dimensional projective variety with $\Pic(V)=1$.
A codimension one holomorphic distribution $\mathcal{D}$ on $V$ of degree $d$ is a global section of $\Omega^1_{V}(d)$.  
\end{Definition}

\begin{thm}
Let $V\subset\mathbb{P}^4$ be a smooth hypersurface, and let
$\mathcal{D}$ be a two-dimensional distribution on $V$ with only isolated singularities.
Then $\mathcal{D}$ is stable in each of the following cases:
\begin{enumerate}
	\item $\deg(V)\neq2$;
	\item $\deg(V)=2$ and $\deg(\mathcal{D})>5$. Alternatively, $\mathcal{D}$ is semi-stable if $\deg(\mathcal{D})\geq5$.
\end{enumerate}
\end{thm}

\noindent Note that $\deg(\mathcal{D})\geq1$, since $H^0\big(V,\Omega^1_{V}(\deg(\mathcal{D}))\big)\neq0$; 
see Subsection \ref{Bottt}.

\begin{proof} Let $V=\left\{f=0\right\}$, where $f$ is a homogeneous polynomial. By Corollary \ref{ton} and Theorem \ref{meme}, 
it suffices to consider the case $\deg(f)=2$; that is, $V=Q_3$ is a quadric.

Using the notation introduced in the proof of Theorem \ref{meme}, where $d\geq2$ and $d-d_1\geq1$, 
Example \ref{qua_n=3} shows that $\wp(2,d)=0$ whenever $d-d_1\leq1$. Hence, we may assume that $d-d_1>1$.
Thus, for Corollary \ref{Exam} $(1)$, we have the following possibilities.
\begin{itemize}
	\item $d<d_1-2$, which is impossible.
	\item $d>d_1+1$, which implies that $2d\geq d+2>d_1+3$. 
	\item $d=d_1$, which implies that $2d>d_1+3$. 
	\item $d=d_1-1$, which implies that $2d=2d_1-2>d_1+3$. 
\end{itemize}
Therefore, $\mu(T_{X})<\mu(T_{\mathcal{D}})$.
\end{proof}

%\newpage
 
%\section{Rascunho}

%\begin{enumerate} 
%	\item Exemplos, variedades de Hopf. 
%	\item Numero de Bruce-Roberts/Tjurina, entre outros indices associados a Flags (para explorar aqui ou noutro projeto).
%\end{enumerate}

%\newpage 


\begin{thebibliography}{99999999}

\bibitem{ams} A. Cavalcante, M. Corr\^ea, S. Marchesi.
{\it On holomorphic distributions on Fano threefolds}
Journal of Pure and Applied Algebra,  Vol. 224, No. 6, (2020).
	
\bibitem{bramafe} J-P Brasselet,  M.  M. Corr\^ ea  ; F. Lourenço: {\it Residues for flags of holomorphic foliations}. 
Advances in Mathematics, Vol. 320, pp 1158-1184 (2017) . 

\bibitem{Ba} P. P. Barbassa: 
{\it Distributions with Unstable Tangent Sheaf on $\mathbb{P}^{3}$}. 
Mediterr. J. Math. (2026) 23:120
May 6, 2026.

\bibitem{BruMe} M. Brunella \and L. G. Mendes: 
{\it Bounding the degree of solutions to Pfaff equations}, 
Publ. Mat. Vol. 44, No. 2, 593-604, (2000).

\bibitem{Bru1} M. Brunella: 
{\it Some remarks on indices of holomorphic vector fields}. Publ. Math. Vol. 41, No. 2, 527-544, (1997).

\bibitem{Can} M. M. Carnicer: 
{\it The Poincar\'e problem in the nondicritical case}. Ann. Math. Vol. 140, No. 2, 289-294, (1994).

\bibitem{CerLN0} D. Cerveau \and A. L. Neto, 
{\it A structural theorem for codimension one foliations on $\mathbb{P}^{n}$
with an application to degree-three foliation}. 
Ann. Sc. Norm. Super. Pisa Cl. Sci. (5), Vol. XII, 1-41, (2013).

\bibitem{Alc} D. Cerveau \and A. Lins Neto: 
{\it Holomorphic foliations in $\mathbb{P}^2$ having an invariant algebraic curve}. 
Ann. Inst. Fourier. Vol. 41, No. 4, 883-903, (1991).

\bibitem{Cerv} D. Cerveau: 
{\it Pinceaux lin\'eaires de feuilletages sur $\mathbb{C}\mathbb{P}(3)$ et conjecture de Brunella}. 
Publ. Mat. 46, no. 2. 441-451 (2002).

\bibitem{CMS} M. Corrêa Jr, L. G. Maza, M. G. Soares,: {\it Algebraic integrability of polyno-
mial differential r-forms}. J. Pure Appl. Algebra 215 (2011), 2290-2294.

\bibitem{Cos} R. Constant da Costa: {\it Splitting aspects of holomorphic distributions with locally free tangent sheaf}.https://arxiv.org/pdf/2405.17415.

\bibitem{OmCoJa1} O. Calvo-Andrade, M. Corr\^ ea \and M. Jardim:  
{\it Codimension One Holomorphic Distributions on the Projective Three-space}. 
International Mathematics Research Notices, Vol. 2020, issue 23, pp. 9011-9074 (2020).

\bibitem{CerLN} D. Cerveau \and A. L. Neto: 
{\it Codimension two holomorphic foliations}. 
J. Differential Geom. 113(3): 385-416 (2019).

\bibitem{SoaCo} M. Corr\^ea \and M. G. Soares: 
{\it Inequalities for characteristic numbers of flags of distributions and foliations}. 
International Journal of Mathematics, Vol. 24, No. 11, 1350093 (2013).

\bibitem{Covi} M. Corrêa, M. Jardim, R. Vidal, {\it On the singular scheme of split foliations} ,
Indiana Univ. Math. J. 64(5) (2015), 1359-1381.

\bibitem{Coma} M. Corrêa Jr, L. G. Maza, M. G. Soares, {\it Hypersurfaces invariant by Pfaff
systems}, Commun. Contemp. Math. 17 (2015), 1450051.

\bibitem{CofeMi} M. Corrêa Jr, A. M. Ferreira, M. Verbitsky: {\it Classification of holomorphic Pfaff systems on Hopf
manifolds}, European Journal of Mathematics (2021) 7:729-740.

\bibitem{Do} I. Dolgachev: 
{\it Weighted projective varieties}. Lecture Notes in Mathematics. Springer Verlag. Vol. 956, (1982).

\bibitem{Este} E. Esteves, S. L. Kleiman, {\it Bounding solutions of Pfaff equations} , Comm.
Algebra 31 (2003), 3771-3793.

\bibitem{FeLo} A. M. Ferreira \and F. Lourenço {\it Baum-Bott residue of flags of holomorphic distributions}, International Journal of Mathematics Vol. 33, No. 13 (2022) 2250087 (17 pages).

\bibitem{Feigin} B. L. Feigin, Characteristic classes of flags of foliations. Funct Anal Its Appl 9, 312317 (1975). 
https://doi.org/10.1007/BF01075877 

\bibitem{Fle} H. Flenner: 
\emph{Divisorenklassengruppen quasihomogener Singularit\"aten}.
J. Reine Angew. Math. 328, 128-160, (1981).

\bibitem{Fi} M. G. Figueira: 
\emph{Extensions and restrictions of holomorphic foliations}.
Annales de la faculté des sciences de Toulouse Volume XXXIII, no 4, 2024
pp. 981-995.

\bibitem{Gal} C. Galindo \and F. Monserrat: 
{\it The Poincar\'e problem, algebraic integrability and dicritical divisors}. 
J. Differ. Equ. 256, 3614-3633, (2014).

\bibitem{HuLe} D. Huybrechts, M. Lehn: 
{\it The Geometry of Moduli Spaces of Sheaves}. Cambridge University Press. 06 July 2010.ISBN:9780511711985.

\bibitem{Ko} J.Kollár et al:
{\it Flips and abundance for algebraic threefolds}.
Astérisque 211 (1993).


\bibitem{Koda} Kodaira, K.:
{\it On the structure of compact complex analytic surfaces II}.Amer. J. Math. 88(3), 682-721 (1966).

\bibitem{HuyLehn} D. Huybrechts and M. Lehn:
{\it The geometry of moduli spaces of sheaves}.2 ed., Cambridge Mathematical Library, Cambridge University Press, 2010.

\bibitem{EvM} E. Materov:
{\it The Bott Formula for Toric Varieties}.
Mosc. Math. J. 2, 161-182, (2002).

\bibitem{mol} R. Mol, {\it Flags of holomorphic foliations } Annals of the Brazilian Academy of Sciences,   83(3): 775-786, 2011.

\bibitem{Mum} D. Mumford:
{\it Projective Invariants of Projective Structures and Applications}.
Proc. Internat. Congr. Mathematicians (1963), 526-530.

\bibitem{LPT} F. Loray, J. V. Pereira and F. Touzet:
{\it Foliations with trivial canonical bundle on Fano 3-folds}.
Math. Nachr. 286, No. 8-9, 921-940, (2013).
 
\bibitem{Pere} J. V. Pereira: 
{\it On the Poincar\'e problem for foliations of general type}. Math. Ann. 323, No. 2, 217-226, (2002).

\bibitem{Poin} H. Poincar\'e: 
{\it Sur l'int\'egration alg\'ebrique des \'equations diff\'erentielles du premier ordre et dupremier degr\'e}. 
Rend. Circ. Mat. Palermo, vol. 5, 161-191, (1891).

\bibitem{Hopf} Hopf, H.: 
{\it Zur Topologie der komplexen Mannigfaltigkeiten}. 
Studies and Essays Presented to R. Courant
on his 60th Birthday, pp. 167-185. Interscience, New York (1948).

\bibitem{Ro}
Rodr\'iguez Pena, Arnulfo Miguel:
{\it On Botts residue formula for toric complete intersections}.
Forum Mathematicum, vol. 38, no. 4, pp. 1079-1098, (2026).

\bibitem{NaTi} A. Napame and C. Tipler:
{\it Toric sheaves, stability and fibrations}. Journal of Pure and Applied Algebra, 228 107494, (2024).

\bibitem{Na} A. Napame:
{\it Stability of equivariant logarithmic tangent sheaves on toric varieties of Picard rank two}.
Annales de la faculté des sciences de Toulouse Volume XXXIII, no 3, 2024
pp. 739-783.

\bibitem{Oko} C. Okonek, M. Schneider, H. Spindler:
{\it Vector bundles on complex projective spaces}.
Progress in Mathematics, vol. 3, Birkhaüser, 1980.

\bibitem{RoFe} A. M. Rodriguez \and F. Lourenco:
{\it Inequalities and enumerative formulas for flags of Pfaff systems}.
Communications in Contemporary Mathematics 2650002, (2026).

\bibitem{Sa} K. Saito: 
{\it On a generalization of De Rham lemma}.
Ann. Inst. Fourier (Grenoble) 26(2), 165-170 (1976).

\bibitem{Sn} D. M. Snow:
{\it Cohomology of twisted holomorphic forms on Grassmann manifolds and quadric hypersurfaces}.
Math. Ann. 276, 159-176 (1986).

\bibitem{So} M. G. Soares: 
{\it The Poincar\'e problem for hypersurfaces
invariant by one-dimensional foliations}. Inventiones Mathematicae. Vol. 128, 495-500, (1997).

\bibitem{Ta} F. Takemoto:
{\it Stable vector bundles on algebraic surfaces}.
Nagoya Math. J. 47 (1972),
29-48.


\end{thebibliography}
\end{document}